\documentclass{amsart}
\usepackage{bbm,amsfonts,latexsym,amsmath,amscd,amssymb,fancyhdr}
\usepackage[all]{xy}

\theoremstyle{plain}
\newtheorem{theorem}{Theorem}
\numberwithin{theorem}{section}

\numberwithin{corollary}{section}

\newtheorem{definition}{Definition}
\numberwithin{definition}{section}

\numberwithin{lemma}{section}

\newtheorem{proposition}{Proposition}
\numberwithin{proposition}{section}

\newtheorem{remark}{Remark}
\numberwithin{remark}{section}

\numberwithin{equation}{section}

\newcommand {\be}{\begin{equation}}
\newcommand {\ee}{\end{equation}}

\newcommand{\h}{\begin{eqnarray*}}
 \newcommand{\e}{\end{eqnarray*}}

\newcommand{\CC}{\mathbf{C}}

\begin{document}
\title[Supersymmetric QFT, Super Loop Spaces and Bismut-Chern Character]{Supersymmetric QFT, Super Loop Spaces and Bismut-Chern Character}
\author{Fei Han}
\address{Department of Mathematics, University of California,
Berkeley, CA, 94720-3840} \email{feihan@math.berkeley.edu}
\address{{\it Current address}: Department of Mathematics, Stanford University, Stanford, CA, 94305-2125} \maketitle

\begin{abstract} In this paper, we give a quantum interpretation of the
Bismut-Chern character form (the loop space lifting of the Chern
character form) as well as the Chern character form associated to a
complex vector bundle with connection over a smooth manifold in the
framework of supersymmetric quantum field theories developed by
Stolz and Teichner \cite{ST07}. We show that the Bismut-Chern
character form comes up via a loop-deloop process when one goes from
$1|1$D theory over a manifold down to a $0|1$D theory over its free
loop space. Based on our quantum interpretation of the Bismut-Chern
character form and Chern character form, we construct Chern
character type maps for SUSY QFTs.
\end{abstract}

\tableofcontents

\section{Introduction}

Roughly speaking, a $d$-dimensional quantum field theory in the
sense of Atiyah and Segal (see \cite{At88}, \cite{Se88}) gives a way
to associate a $(d-1)$-dimensional manifold a Hilbert space and to a
bordism between such manifolds a trace class operator between the
corresponding Hilbert spaces, such that gluing of bordisms
corresponds to composing operators. A field theory based on a
manifold $M$ is as before except that the $(d-1)$-manifolds and the
bordisms between them come equipped with maps to $M$. Segal
\cite{Se04} suggests that the conformal field theories (CFT) over a
manifold $M$ might be able to provide cocycles for the elliptic
cohomology of $M$. Stolz and Teichner \cite{ST04} have developed
Segal's idea by adding world-sheet locality and supersymmetry (SUSY)
into the picture. They conjecture that the space of certain 2D SUSY
QFTs gives the spectrum TMF of topological modular forms \cite{Ho}
and for a smooth manifold $M$, the space of all 2D SUSY QFTs over
$M$ moduled out proper concordance relations gives
$\mathrm{Ell}^n(M)$, the elliptic cohomology of $M$. See also
\cite{BDR02}, \cite{HK} and \cite{DLM05} for other contributions on
the aim to understand elliptic cohomology in geometric ways.

The K-theory can be considered as a ``case study" for elliptic
cohomology. Stolz and Teichner \cite{ST04} have shown that
1-dimensional SUSY QFTs indeed give the spectrum for K(or
KO)-theory. To geometrically understand their 1-dimensional theory,
one has to study SUSY $1$D field theories over parametrizing
manifolds. Along this direction, Florin Dumitrescu \cite{D06} has
developed the theory of super parallel transport. More precisely,
given a vector bundle $E$ with a connection $\nabla^E$ over a
manifold $M$, let $S$ be a super manifold and $c: S\times
\mathbf{R}^{1|1}\rightarrow M$ be a $S$-family super path in $M$, he
is able to construct a bundle map $SP(c,t):c^*_{0,0}E\rightarrow
c^*_{t, \theta}E$, which has similar nice properties as the ordinary
parallel transport. In chapter 3, we will briefly introduce (and a
little bit modify) his construction to obtain a functor from the
category of smooth complex vector bundles over $M$ with connections
to the category of SUSY 1D QFTs over $M$. Florin's construction
provides examples of $1|1$D SUSY QFTs.

Hohnhold, Kreck, Stolz and Teichner have also studied $0|1$D SUSY
QFTs over $M$ and related them to the theory of differential forms
and de Rham cohomology \cite{HKST06}. More precisely, they show that
the set of certain 0D SUSY QFT's is canonically mapped to the set of
closed differential forms by a bijection. Concordance between the
field theories corresponds to the differential forms being
cohomologous to each other.

Stolz and Teichner conjecture that there should be a beautiful
quantum interpretation of the Chern character in terms of a map from
$1|1$D QFTs over $M$ to $0|1$D QFTs over the free loop space $LM$,
given by crossing with the standard circle.

We confirm their conjecture in this paper. Starting from a vector
bundle $E$ with connection $\nabla^E$ over $M$, applying
Dumitrescu's construction of super parallel transport \cite{D06},
one obtains a SUSY $1|1$D QFT over $M$. Choosing a special
$S$-family super path with $S$ the super loop space and $c$ the
super evaluation map, we obtain a special super parallel transport.
We then use this super parallel transport to construct a
differential form on $LM$ via some loop-deloop process.  It turns
out that this differential form is just the Bismut-Chern character
form \cite{B85} over $LM$. Note that Bismut obtained it by extending
the ideas of Witten and Atiyah \cite{At84} of interpreting the index
of the Dirac operator on the spin complex of a spin manifold as a
paring of $1\in \Omega(LM)$ with certain equivariant current $\mu_D$
on the loop space to the situation of twisted spin complex. The
restriction of the Bismut-Chern character form to $M$, the $S^1$
fixed point set on $LM$, is the ordinary Chern character form
associated to $(E, \nabla^E)$. See also \cite{GJP91} for a
representation of the Bismut-Chern character form in the cyclic bar
complex.

Our construction of the Bismut-Chern character form shows that it
actually represents a loop-deloop process in the framework of
supersymmetric quantum field theories (see Theorem 4.6 and Section
4.4 in this paper for details). It provides us a new way to
understand the Bismut-Chern character form over $LM$ as well as the
Chern character form over $M$. The character form can be considered
as a phenomena when one goes from $1|1$D field theories to $0|1$D
field theories. We hope to understand this kind of phenomena when
the dimensions of the field theories are higher. What we did
actually gives the Bismut-Chern character form as well as the Chern
character form a quantum interpretation (Theorem 4.6 and Theorem
4.7). Based on this quantum interpretation and the 0-dimensional
theory in \cite{HKST06}, we are able to construct Chern character
type maps for SUSY QFTs.

The paper is organized as follows. In section 2, we supply necessary
preliminary knowledge on super geometry. In section 3, we briefly
introduce SUSY QFTs developed by Stolz and Teichner as well as
Dumitrescu's super parallel transport and its slight modification.
We at last construct the Bismut-Chern character as well the Chern
character and then construct Chern character type maps in the
framework of SUSY QFTs in section 4.

\section{Preliminaries on Super Geometry} In this section, we
briefly survey the theory of super manifolds. We refer interested
readers to the standard references \cite{DM99} and \cite{Va04} for
details.

A {\it supermanifold} $M$ of dimension $(m|n)$ is a pair $(|M|,
\mathcal{O}_M)$, consisting of \newline \noindent 1) an $m$
dimensional smooth manifold $|M|$, the so called underlying
manifold, or the reduced manifold; \newline \noindent 2) a sheaf
of $\mathbf{Z}_2$ graded commutative algebras on $|M|$, the
``structure sheaf" $\mathcal{O}_M$, such that for any open set $U
\subset |M|$,
$$\mathcal{O}_M(U)\cong C^\infty(U)\otimes \Lambda [\theta_1,
\theta_2, \cdots, \theta_n].$$

We refer to global sections of the structure sheaf of $M$ as
elements in $C^\infty(M):=\Gamma(\mathcal{O}_M).$

A {\it morphism between two super manifolds} $M=(|M|,
\mathcal{O}_M), N=(|N|, \mathcal{O}_N)$ is a pair $(f,
f^\natural)$ such that
$$f:|M|\rightarrow |N|$$ is a smooth map between the underlying
manifolds and
$$f^\natural: \mathcal{O}_N \rightarrow f_*\mathcal{O}_M$$ is a
morphism of sheaves.

Let $\mathbf{SM}$ denote the category of supermanifolds and
$\mathbf{SM}(M,N)$ denote the morphisms between $M$ and $N$. A basic
fact in supergeometry is that maps between two supermanifolds are
uniquely determined by the map induced on global sections (see
\cite{Kos77}, page 208). So we quite often write a map between
sheaves as just the map induced on their global sections.

Let $E$ be a vector bundle over an ordinary manifold $M$. One can
canonically construct a super manifold $\Pi E=(M, \mathcal{O}_{\Pi
E})$ where $\mathcal{O}_{\Pi E}$ is the sheaf of sections of
$\Lambda E^*$. This construction provides an important source of
supermanifolds. This construction also defines a functor:
$$S: \mathbf{VB}\rightarrow \mathbf{SM}: E\mapsto \Pi E.$$ This
functor actually induces a bijection on the isomorphism classes of
objects however it does not give an equivalence of categories
because there are more more morphisms in $\mathbf{SM}$ then in
$\mathbf{VB}$ (cf. \cite{D06}).

Since sheaves are generally hard to work with, one often thinks of
super manifolds in terms of their ``$S$-points". To be precise,
instead of $M$ itself, one considers the morphisms sets
$\mathbf{SM}(S, M)$, where $S$ varies over all super manifolds
$S$. One can think of an $S$-point as a family of points of $M$
parametrized by $S$. It's not hard to see that $M$ determines a
contravariant functor:
$$\mathbf{SM}\rightarrow \mathbf{Sets}: S \mapsto
\mathbf{SM}(S,M).$$ It's called the {\it functor of points} of $M$.
A map $f:M\rightarrow N$ of super manifolds determines a natural
transformation $\mathbf{SM}(\cdot, M)\rightarrow \mathbf{SM}(\cdot,
N)$. The Yoneda's lemma (cf. \cite{N96}) tells us that
$$ \mathbf{SM} \rightarrow \mathcal{F}(\mathbf{SM}^{op}, \mathbf{Sets}),\  M \mapsto \mathbf{SM}(\cdot,
M)$$ is an embedding of categories. This means that to give a map
$M\rightarrow N$ amounts to give a natural transformation of
functors $\mathbf{SM}(-, M)\rightarrow \mathbf{SM}(-, N)$. From
this point of view, one can think of supermanifold $M$ as a
representable functor $\mathbf{SM}\rightarrow \mathbf{Sets}$. Such
a functor determines $M$ uniquely up to isomorphism. The {\it
product} $M\times N$ can be interpreted as the supermanifold
representing the functor $\mathbf{SM}(\cdot, M)\times
\mathbf{SM}(\cdot, N).$

An arbitrary contravariant functor $\mathbf{SM}\rightarrow
\mathbf{Sets}$ will be called a {\it generalized supermanifold}.
Therefore Yoneda's lemma says $\mathbf{SM}$ embeds faithfully into
the category $\mathbf{GSM}$ of generalized supermanifolds.

Given two supermanifolds $M, N$, consider the generalized
supermanifold
$$ \underline{\mathbf{SM}}(M, N): \mathbf{SM}\rightarrow
\mathbf{Sets}: S\rightarrow \mathbf{SM}(S\times M, N).$$ If
$\underline{\mathbf{SM}}(M, N)$ is an ordinary supermanifold, then
by definition we have the following adjunction formula \be
\mathbf{SM}(S,\underline{\mathbf{SM}}(M, N))\cong
\mathbf{SM}(S\times M, N).\ee

Let $M$ and $N$ be two supermanifolds, define the following
``evaluation" map \be ev:\underline{\mathbf{SM}}(M, N)\times M
\rightarrow N\ee via it's $S$-points. That means that for any
supermanifold $S$, define  \be ev_S:\mathbf{SM}(S\times M, N)\times
\mathbf{SM}(S, M)\rightarrow \mathbf{SM}(S, N), \ee \be
ev_S(\varphi, m)=\varphi \circ (1\times m)\circ \Delta \in
\mathbf{SM}(S, N).\ee

\begin{remark}Although supermanifolds have rigorous mathematical
foundations now, in practice, physicists often freely and correctly
treat supermanifolds as ordinary manifolds with `` even and odd
variables" without always coming back to the rigorous mathematical
definitions.
\end{remark}

On a supermanifold, one can define {\it tangent sheaf} and {\it
tangent vectors} on it. The tangent sheaf $\mathcal{TM}$ is defined
as the sheaf of graded derivations of $\mathcal{O}_M$, i.e. for
$U\subseteq |M|$
$$\mathcal{TM}(U)=\{X:\mathcal{O}(U)\rightarrow \mathcal{O}(U)\
\mathrm{linear}\,:X(fg)=X(f)g+(-1)^{p(X)p(f)}fX(g)\}.$$ Here
$p(X)=0$ or $1$ according to whether $X$ is even, respectively odd
vector field on $U$, and similarly $p(f)=0 \ \mathrm{or}\  1$, for
$f$ even, respectively odd, function on $M$. $\mathcal{TM}$ is then
a locally free $\mathcal{O}_M$-module of rank $(p,q)$ the dimension
of the supermanifold $M$. Sections of $\mathcal{TM}$ are the vector
fields on $M$. For $X$ and $Y$ vector fields on $M$, define their
{\it Lie bracket} $[X,Y]$ by
$$[X,Y](f)=X(Y(f))-(-1)^{p(X)p(Y)}Y(X(f))$$
for $f\in C^\infty(M)=\mathcal{O}_M(|M|).$

For example, on $\mathbf{R}^{1|1}$, there are two canonical vector
fields $D=\frac{\partial}{\partial
\theta}+\theta\frac{\partial}{\partial t}$ and
$Q=\frac{\partial}{\partial \theta}-\theta\frac{\partial}{\partial
t}.$ It's not hard to verify that \be
D^2=\frac{1}{2}[D,D]=\frac{\partial}{\partial t}, \
Q^2=\frac{1}{2}[Q,Q]=-\frac{\partial}{\partial t}, \ [D, Q]=0. \ee

Dually one can also define {\it cotangent sheaf} and {\it exterior
cotangent sheaf} on a supermanifold. Define the cotangent sheaf
$\Omega^1(M)$ to be the dual of the tangent sheaf $\mathcal{TM}$.
Sections of $\Omega^1(M)$ are called differential one forms. Let
$\langle\ ,\  \rangle: \mathcal{TM}\times \Omega^1(M)\rightarrow
\mathcal{O}_M$ denote the duality pairing between vector fields and
1-forms. Define the {\it exterior derivative}
$d:\mathcal{O}_M\rightarrow \Omega^1(M)$ by
$$\langle X, df\rangle=X(f), \ \mathrm{for}\  X\in \mathcal{TM}, f\in
\mathcal{O}_M.$$ Let $\Omega^*(M)=\Lambda^*\Omega^1(M)$ be the
exterior cotangent sheaf on $M$, whose sections are called
differential forms on $M$. $d$ extends uniquely to a degree one
derivation $d:\Omega^*(M)\rightarrow\Omega^*(M)$ by requiring that
$$ d^2=0,$$
$$d(\alpha \wedge \beta)=d\alpha\wedge \beta+(-1)^p\alpha\wedge
d\beta, \ \mathrm{for} \ \alpha\in\Omega^p(M).$$

The readers can consult \cite{DM99} for more details.

One has the following proposition (cf. \cite{D06}, \cite{HKST06} for
a proof)
\begin{proposition} Let $M$ be an ordinary manifold. Then we can
identify \be \underline{\mathbf{SM}}(\mathbf{R}^{0|1}, M)\cong \Pi
TM.\ee
\end{proposition}
\begin{proof}We want to show that we have isomorphisms
$$ \Psi_S: \mathbf{SM}(S\times \mathbf{R}^{0|1}, M)\to
\mathbf{SM}(S, \Pi TM),$$ natural in $S$, where $S$ is an arbitrary
supermanifold. The left hand side is the set of grading preserving
maps of $\mathbf{Z}_2$-algebras
$$\varphi: C^\infty(M)\to C^\infty(S\times
\mathbf{R}^{0|1})=C^\infty(S)[\theta].$$ If we write
$\varphi(f)=\varphi_1(f)+\theta\varphi_2(f)$, for $f\in
C^\infty(M)$, then the fact that $\varphi(fg)=\varphi(f)\varphi(g)$
is equivalent to the following conditions:
$$\varphi_1(fg)=\varphi_1(f)\varphi_1(g), \ \ \ \varphi_2(fg)=\varphi_2(f)\varphi_1(g)+(-1)^{p(f)}\varphi_1(f)\varphi_2(g).  $$
The first condition is equivalent to $\varphi_1=a^\natural$, for
some $a:S \to M$. The second condition tells us that $\varphi_2$ is
an odd tangent vector at $a\in M(S)$, i.e. $\varphi_2=X_a\in TM_a$.
Therefore the left hand side is
$$\mathbf{SM}(S\times \mathbf{R}^{0|1}, M)=\{\mathrm{pairs}\  (a,
X_a)|a\in M(S), X_a\in TM_a, X_a \ \mathrm{odd}\}.$$

The right hand side $\mathbf{SM}(S, \Pi TM)$ is the set of
$\mathbf{Z}_2$-graded algebra maps $\Omega^*(M)\to C^\infty(S).$
Such maps are determined by their restriction to 0-forms and
1-forms. Define then $\Psi_S(a, X_a)$ to be the map $S\to \Pi TM$
determined by defining it on functions $f\in C^\infty(M)$ by
$a^\natural(f)\in C^\infty(S)$ and on forms $\omega$ by $\langle
X_a, \omega \rangle \in C^\infty(S)$. One can easily check that
$\Psi_S$ is well-defined, bijective and natural in $S$.
\end{proof}

\section{SUSY QFTs} In this section, we give a very rough introduction to formulations of SUSY QFT's over manifolds developed by Stolz and Teichner
(\cite{HKST06}, \cite{HST07}, \cite{ST06}, \cite{ST07}) as well as
Dumitrescu's construction of super parallel transport \cite{D06} and
our slight modifications of his construction. Such formulations of
supersymmetric quantum field theories are cocycles for generalized
cohomology theories via SUSY QFTs. The purpose of this section is to
put our construction of Bismut-Chern character in the next section
into this framework of SUSY QFT's.

Quantum field theories in the sense of Atiyah-Segal is a functor
from a suitable bordism category to a category of locally convex
topological vector spaces satisfying certain axioms. Following Stolz
and Teichner, we will briefly first introduce the rendition of
Atiyah-Segal's quantum field theories and further enrich it by
adding smoothness and supersymmetry. For complete details of all
those categories and functors, see \cite{ST07}.

\subsection{Preliminary Definition of QFTs}

Let's first introduce relevant categories.
\begin{definition} {\bf (The Riemannian spin bordism
category RB$^d$)}The objects and morphisms in the category RB$^d$
are as follows:\newline
 \noindent {\bf objects} are quadruples (U, Y, U$^{-}$, U$^{+}$),
 where
 \begin{itemize}
\item[]
 U is a Riemaannnian spin manifold of dimension d (typically not
 closed);\\[2mm]
 Y is a closed codimension 1 smooth submanifold of $U$;\\[2mm]
 U$^{\pm}$ are disjoint open subsets of $U\setminus Y$ whose union is $U\setminus Y$. Y is
 r required to contain in the closure of both U$^{\pm}$ (this
 ensures that U$^{\pm}$ are collars of Y).
\end{itemize}
We will often suppress (U, Y, U$^-$, U$^+$) in the notation and just
write Y instead of (U, Y, U$^-$, U$^+$).

\noindent {\bf morphisms} from $Y_1$ to $Y_2$ are equivalence
classes of Riemannian spin bordisms from $Y_1$ to $Y_2$; here
Riemannian spin bordism is a triple $(\Sigma, \iota_1, \iota_2)$,
where
 \begin{itemize}
\item[]
 $\Sigma$ is a Riemannian spin manifold of dimension d (not necessarily closed),
 and\\[2mm]
 $\iota_1:V_1\hookrightarrow \Sigma$ and $\iota_2:V_2\hookrightarrow
 \Sigma$ are isometric spin embeddings, where $V_k\subset U_k$ for
 k=1,2 is some open neighborhood of $Y_k\subset U_k$.
\end{itemize}
We define $V_k^{\pm}:=U_k^{\pm}\cap V_k$ and require that
\begin{itemize}
\item[]
 $\bullet$ the sets $\iota_1(V_1^+\cup Y_1)$ and $\iota_1(V_2^-\cup Y_2)$ are disjoint and\\[2mm]
 $\bullet$ $\Sigma \setminus (\iota_1(V_1^+)\cup \iota_2(V_2^-))$ is
 compact.
\end{itemize}
Note that $\Sigma \setminus (\iota_1(V_1^+)\cup \iota_2(V_2^-))$ is
a compact manifold with boundary $\iota_1(Y_1)\coprod \iota_2(Y_2)$;
i.e., it's a bordism between $Y_1$ and $Y_2$ in the usual sense. Now
suppose that $(\sum, \iota_1, \iota_2)$ and $(\Sigma', \iota_1',
\iota_2')$ are two Riemaniann spin bordisms from $Y_1$ to $Y_2$ with
$V_1\subset V_1', V_2\subset V_2'$ and that there is a spin isometry
$F$ that makes the following diagram commutative

\begin{displaymath}
\xymatrix{
  V_2 \ar[r]^{\iota_2} \ar[d]_{i} &\Sigma  \ar[d]_{F}^{\cong}  &   V_1 \ar[l]_{\iota_1} \ar[d]^{i}   \\
    V_2' \ar[r]^{\iota_2'}                          &\Sigma'  &   V_1' \ar[l]_{\iota_1'}                }
\end{displaymath} Then we declare that $(\Sigma, \iota_1, \iota_2)$
as equivalent to $(\Sigma', \iota_1', \iota_2')$. A morphism from
$Y_1$ to $Y_2$ is an equivalence class of Riemannian spin bordisms
with this equivalence relation.

\noindent {\bf composition} of Riemannian spin bordisms is given by
gluing; more precisely, let $(\Sigma', \iota_1', \iota_2')$ be a
Riemannian spin bordism from $Y_1$ to $Y_2$ and $(\Sigma, \iota_2,
\iota_3)$ a Riemannian spin bordism from $Y_2$ to $Y_3$. Without
loss of generality, we can assume that the domains of the isometries
$\iota_2'$ and $\iota_2$ agree; suppose that $V_2\subset U_2$ is
this common domain. Then identifying $\iota_2'(V_2)\subset \Sigma'$
with $\iota_2(V_2)\subset \Sigma$ via the isometry
$\iota_2\circ(\iota_2')^{-1}$ gives the Riemaniann spin manifold
$\Sigma'':=\Sigma\cup_{V_2}\Sigma'$.
\end{definition}

There are some additional structures on the bordism category RB$^d$.
They are
\begin{itemize}
\item[]
{\bf symmetric monoidal structure}: disjoint union gives RB$^d$ the
structure of a symmetric monoidal category; the unit
object is given by the empty $(d-1)$-manifold.\\[2mm]
{\bf the anti-involution ${}^\vee$} On objects the anti-involution
${}^\vee$ is defined by interchanging $U^+$ and $U^-$ (which can be
thought of as flipping the orientation of the normal bundle to $Y$
in $U$). If $(\Sigma, \iota_1, \iota_2)$ is a Riemannian bordism
from $Y_1$ to $Y_2$, then $(\Sigma, \iota_1, \iota_2)^\vee=(\Sigma,
\iota_2, \iota_1)$ is a Riemannian bordism from $Y_2^\vee$ to
$Y_1^\vee$.\\[2mm]
{\bf the involution {}$^-$}: Replacing the spin structure on the
bicollars $U$ as well as the bordism $\Sigma$ by their opposite
defines an involution {}$^-$: RB$^d$ $\rightarrow$ RB$^d$.
\end{itemize}

The other categories involved in QFT are TV and TV$^{\pm}$.

\begin{definition}{\bf (the category TV)}\newline
{\bf objects} are $\mathbf{Z}/2$-graded locally convex vector
spaces;\newline {\bf morphisms} are graded preserving continuous
linear maps.
\end{definition}

\begin{definition}{\bf (the category TV$^{\pm}$)}\newline
{\bf objects}are triples $V=(V^+, V^-, \mu_V)$, where $V^{\pm}$ are
locally convex vector spaces, and $\mu_V:V^-\otimes V^+\rightarrow
\mathbf{C}$ is a continuous linear map. Here and in the following,
$\otimes$ is the projective tensor product.
\newline {\bf morphisms} from $V=(V^+, V^-, \mu_V)$ to $W=(W^+, W^-,
\mu_W)$ are pairs $T=(T^+:V^+\rightarrow W^+, T^-:W^-\rightarrow
V^-)$ of continuous linear maps, which are dual to each other in the
sense that
$$ \mu_V(T^-w^-\otimes v^+)=\mu_W(w^-\otimes T^+v^+), \forall v^+\in
V^+, w^-\in W^-.$$ \noindent {\bf composition} Let $S=(S^+, S^-):
V_1\rightarrow V_2$ and $T=(T^+, T^-): V_2\rightarrow V_3$ be
morphisms in TV. Then their composition $T\circ S: V_1\rightarrow
V_3$ is given by
$$ T\circ S:=(T^+\circ S^+, S^-\circ T^-).$$
\end{definition}

There are also some additional structures on the category
TV$^{\pm}$. They are
\begin{itemize}
\item[]
{\bf symmetric monoidal structure}: The tensor product of two objects $V=(V^+, V^-, \mu_V)$ and
$W=(W^+, W^-, \mu_W)$ is defined as
$$ V\otimes W:=(V^+\otimes W^+, V^-\otimes W^-, \mu_{V\otimes W}),$$
where $V^{\pm}\otimes W^{\pm}$ is the projective tensor product and
$\mu_{V\otimes W}$ is given by the composition of the usual graded
symmetry isomorphism
$$ (V^-\otimes W^-)\otimes (V^+\otimes W^+)\cong V^-\otimes
V^+\otimes W^-\otimes W^+$$ and the linear map
$$\xymatrix@C=0.5cm{
 V^-\otimes V^+\otimes W^-\otimes
W^+ \ar[rr]^{\ \ \ \ \ \ \  \ \ \ \ \ \mu_V\otimes \mu_W} &&
\mathbf{C}\times \mathbf{C}}=\mathbf{C}.$$ On morphisms, we define
$(T^+, T^-)\otimes (S^+, S^-):=(T^+\otimes S^+,
T^-\otimes S^-).$\\[2mm]
{\bf the anti-involution ${}^\vee$} On objects On objects, it is given by $(V^+, V^-, \mu_V)^\vee:=(V^-, V^+, \mu_V^\vee),$
where $\mu_V^\vee$ is the composition
$$\xymatrix@C=0.5cm{
  V^+\otimes V^- \ar[rr]^{\cong} && V^-\otimes V^+ \ar[rr]^{\ \ \ \ \ \ \mu_V} && \mathbf{C}}$$ of the graded symmetry isomorphism and $\mu_V$. On
morphisms the anti-involution is given by $(T^+, T^-)^\vee:=(T^-, T^+)$. We note that $V^\vee\otimes W^\vee=(V\otimes W)^\vee$.\\[2mm]
{\bf the involution {}$^-$}: If $V=(V^+, V^-, \mu_V)$ is an object
of TV$^{\pm}$, then $\overline{V}$ is given by complex conjugate
vector spaces $\overline{V}^+, \overline{V}^-$ and the paring
$$\overline{V}^-\otimes \overline{V}^+=\overline{V^-\otimes V^+}\xymatrix@C=0.5cm{
 \ar[rr]^{\overline{\mu}_V} && \overline{\mathbf{C}}\cong
 \mathbf{C}}.$$ On morphisms, it is given by $\overline{(T^+, T^-)}=(\overline{T}^+,
 \overline{T}^-)$, where (as for $\mu_V$) $\overline{T}^{\pm}$ is
 the same map as $T^{\pm}$, but regarded as a complex linear map
 between the complex conjugate vector spaces.
\end{itemize}

There is also a functor TV$\rightarrow$TV$^\pm$; on objects, it
sends a locally convex vector space to $(V, V', \mu)$, where $\mu:
V'\otimes V\rightarrow \mathbf{C}$ is the natural paring, $V'$ is
the continuous dual of $V$. On morphisms, it sends a linear map $T:
V\rightarrow W$ to the pair $(T, T')$, where $T':W'\rightarrow V'$
is the continuous dual to $T$. This is not a {\it monoidal} functor
due to the incompatibility of $\otimes$ and $'$ mentioned earlier.
However one obtains monoidal functor if restricting to {\it finite
dimensional} vector spaces. Using this functor, we can interpret
finite dimensional vector spaces and linear maps as objects resp.
morphisms in TV$^{\pm}$.

With the above preparations, we can give the preliminary definition
of QFT.
\begin{definition}A quantum field theory of dimensional d is a
symmetric monoidal functor
$$F^d: \mathrm{RB}^d\rightarrow \mathrm{TV}^{\pm}$$ which is
compatible with the involution ${}^-$ and the anti-involution
${}^\vee$.
\end{definition}

\subsection{QFTs as Smooth Functors}

Now let's enrich the above definition of QFT by adding {\bf
smoothness} into the picture.

\begin{definition}{\bf (Smooth categories and functors)} A smooth
category $\mathfrak{C}$ is a functor
$$\mathfrak{C}:\mathrm{MAN}^{op}\rightarrow \mathrm{CAT}$$ from the
category of smooth manifolds to the category of categories. If $S$
is a manifold, we will write $\mathfrak{C}_S$ for the category
$\mathfrak{C}(S)$; if $f: S'\rightarrow S$ is a smooth map, we will
write $f^*: \mathfrak{C}_S \rightarrow \mathfrak{C}_{S'}$ for the
corresponding functor $\mathfrak{C}(f)$.

If $\mathfrak{C}$ and $\mathfrak{D}$ are smooth categories, a smooth
functor $\mathfrak{F}$ from $\mathfrak{C}$ to $\mathfrak{D}$ is a
natural transformation. For a manifold $S$, we will write
$\mathfrak{F}_S: \mathfrak{C}_S\rightarrow \mathfrak{D}_S$ for the
corresponding functor.
\end{definition}

In the following, we proceed to define the smooth version of the
above categories and functors involved in the above definition of
QFT. Let's first introduce the notion of quasi bundles (the purpose
of this concept is to guarantee $(fh)^*=g^*f^*$ for pull backs).

\begin{definition} A quasi bundle over a smooth manifold $S$ is a
pair $(h, V)$, where $h:S\rightarrow T$ is a smooth map and
$V\rightarrow T$ is a smooth , locally trivial bundle over $T$. If
$(h', V')$ is another quasi bundle over $S$, a map from $(h,V)$ to
$(h', V')$ is a smooth bundle map $F:h^*V \to (h')^*V'$.
\end{definition}

A smooth map $f:S'\to S$ induces a contravariant functor
$$f^*: \mathrm{QBUN}(S)\to \mathrm{QBUN}(S')$$ from the category of
quasi bundles over $S$ to those over $S'$. It is defined by $$f^*(h,
V)=(f\circ h, V)$$ on objects; on morphisms it is given by the usual
pullback via $f$. In particular, the functor $(fg)^*$ is equal to
$g^*\circ f^*$. Note that the category of quasi bundles over $X$ is
equivalent to the category of bundles over $X$, sending $(h, V)$ to
$h^*V$ provides the equivalence.

The $S$-family version of the RB$^d$ is the following.

\begin{definition}{\bf (The category $\mathrm{RB}^d_S$)} Let $S$ be a
smooth manifold. We want to define the category RB$^d_S$ of
$S$-families of Riemannian bordisms of dimension d in such a way
that RB$^d_{pt}$ agree with the category RB$^d$ of the above
definition 3.1. In that definition, objects and morphisms were
defined in terms of the category Riem$^d$ whose objects are
Riemannian spin manifolds of dimension d and morphisms are isometric
spin embeddings. Here we replace Riem$^d$ by Riem$^d_S$, whose
objects (resp. morphisms) are S-families of objects (resp.
morphisms) of Riem$^d$. More precisely,
\begin{itemize}
\item[]
$\bullet$ an object of Riem$^d_S$ is a smooth quasi bundle $U\to S$ with d-dimensional fibers, euqipped
with a fiberwise Riemannian metric and spin structure (i.e. a spin structure on the vertical tangent
bundle).\\[2mm]
$\bullet$ A morphism from $U$ to $U'$ is a smooth quasi bundle map
$f: U\to U'$ preserving the fiberwise Riemannian metric and spin
structure.
\end{itemize}
\end{definition}

A smooth map $f:S'\to S$ induces a pullback functor $f^*:
\mathrm{RB}^d_S\to \mathrm{RB}^d_{S'}$, such that $(fg)^*=g^*f^*$.
The fiberwise disjoint union gives RB$^*$ the structure of a
symmetric monoidal category; the involution ${}^-$, the
anti-involution ${}^\vee$ and the adjoint transformation generalize
from RB$^d$ to RB$^d_S$.

\begin{definition}{\bf(The smooth category $\mathfrak{RB}^d$)} The
smooth category $\mathfrak{RB}^d$ is the functor
$$\mathfrak{RB}^d: \mathrm{MAN}^{op}\rightarrow \mathrm{CAT}$$ which sends a smooth
manifold $S$ to the category $RB_S$ and a smooth map $f:S\to S'$ to
the pullback functor $f^*: RB^d_S\to RB^d_{S'}$.
\end{definition}

The $S$-family version of the TV$^{\pm}$ is the following.
\begin{definition} {\bf (The category $\mathrm{TV}^{\pm}_S$)} Let $S$
be a smooth manifold. We want to define the category TV$^{\pm}_S$ of
S-families of (pair) of locally convex vector space in such a way
that $TV^{\pm}_{pt}$ agree with the category $TV^{\pm}$ of
definition 3.3. Just replace locally convex vector spaces by quasi
bundles of locally convex vector spaces over S and continuous linear
maps by smooth bundle maps. In particular, objects are triples
$V=(V^+, V^-, \mu_V)$, where $V^{\pm}$ are quasi bundles of locally
convex vector spaces over $S$, and $\mu_V: V^-\otimes V^+\to \CC_S$
is a smooth bundle map. Here the tensor product of quasi bundles
over $S$ is given by
$$(h: S\to T, V)\otimes (h': S\to T', V'):=(h\times h':S\to T\times
T', p_1^*V\otimes p_2^*V'),$$ where $p_1^*V\otimes p_2^*V'\to
T\otimes T'$ is a fiberwise (projective) tensor product and $p_1$
(res. $p_2$) is the projection onto the first (res. second) factor.
Moreover, $\CC_S$ is the quasi bundle $(p:S\to pt, pt\times \CC)$;
replacing $\CC$ by any locally convex vector space $V$ gives a quasi
bundle $V_S$ over $S$. We note that $C_S$ is the unit for the tensor
product of quasi bundles over $S$.
\end{definition}

A smooth map $f:S'\to S$ induces a functor
$$ f^*:\mathrm{TV}^{\pm}_S\to \mathrm{TV}^{\pm}_{S'}$$ via pullback
of quasi bundle. If $g:S''\to S'$ is a smooth map, then the functor
$(gf)^*$ is equal to $g^*f^*$. The additional structures for the
category TV$^{\pm}$ (the symmetric monoidal structure, the
(anti-)involution and the adjunction transformation) generate in a
straightforward way to the category TV$^{\pm}_S$. These structures
are compatible with pullback functor $f^*$.

\begin{definition}{\bf (The smooth category $\mathfrak{TV}^{\pm}$)}
The smooth category $\mathfrak{TV}^{\pm}$ is the functor
$$\mathfrak{TV}^{\pm}: \mathrm{MAN}^{op}\to \mathrm{CAT}$$ which
sends a smooth manifold $S$ to the category TV$^{\pm}_S$ and a
smooth map $f:S'\to S$ to the pullback functor
$f^*:\mathrm{TV}^{\pm}_S\to \mathrm{TV}^{\pm}_{S'}$.
\end{definition}

\begin{definition} A quantum field theory of dimension d is a
natural transformation
$$ \mathfrak{F}^d: \mathfrak{RB}^d\to \mathfrak{TV}^{\pm},$$ such that for any $S\in \mathrm{MAN}$, $\mathfrak{F}^d_S$ is
a symmetric monoidal functor, which is compatible with the
involution ${}^-$ and the anti-involution ${}^{\vee}$.
\end{definition}

There is also a relative version of the above story. Let $M$ be a
smooth manifold. Define $\mathrm{RB}^d(M)$ to be the category with
objects $(U, Y, U^+, U^-, \alpha:U\to M)$, where $\alpha:U\to M$ is
a smooth map; morphisms $(\Sigma, \iota_1, \iota_2, \beta: \Sigma\to
M)$, where $\beta:\Sigma\to M$ is a smooth map. Those relative maps
satisfy natural commutativity conditions when we do bordisms.
Similarly, one can define family versions $\mathrm{RB}_S^d(M)$ and
$\mathfrak{RB}^d(M)$.  We define a {\bf quantum field theory of
dimension $d$ over $M$} to be a natural transformation
$$ \mathfrak{F}^d(M): \mathfrak{RB}^d(M)\to \mathfrak{TV}^{\pm},$$
 such that for any $S\in \mathrm{MAN}$, $\mathfrak{F}^d_S(M)$ is
a symmetric monoidal functor, which is compatible with the
involution ${}^-$ and the anti-involution ${}^{\vee}$.

\subsection{SUSY QFTs and Examples} Our next step is to enrich our definitions of QFT
described above by adding supersymmetry into the picture.

\subsubsection{Definitions of SUSY QFTs} In the following sections, we will always work with
{\bf complex super manifolds} or {\bf cs-manifolds} (cf.
\cite{DM99}). A {\bf cs-manifold of dimension $n|m$} is a
topological space $M_{red}$ together with a sheaf $\mathcal{O}_M$ of
graded commutative algebras over complex numbers, which is locally
isomorphic to $\mathbf{R}^{n|m}_{cs}:=(\mathbf{R}^{n|m},
\mathcal{O}^{n|m}\otimes \mathbf{C})$. If $M=(M_{red}, \mathcal{O})$
is a cs-manifold, we denote by $\overline{M}:=(M_{red},
\overline{\mathcal{O}})$ the complex conjugate cs-manifold, where
$\overline{\mathcal{O}}(U)=\overline{\mathcal{O}(U)}$. A super
manifold of dimension $n|m$ leads to a cs-manifold by complexifying
its structure sheaf.

To give the definition of SUSY QFT, we will also have to define the
so called {\bf super Riemannian structures} on cs-manifolds
dimension $d|1$. However we will only describe super Riemannian
structures on cs-manifolds of dimension $1|1$ for the purpose for
this paper. For general definitions and their physics motivations,
the readers are referred to \cite{ST07}.

\begin{definition} Let $M$ be a cs-manifold of dimension $1|1$. A
super Riemannian structure on $M$ is given by a collection of pair
$(U_i, D_i)$ indexed by some set $I$, where \newline 1) the $U_i$'s
are open subsets of $M_{red}$ whose union is all of $M_{red}$.
\newline 2) the $D_i$'s are sections of the tangent sheaf $TM$
restricted to $U$ satisfying
\begin{itemize}
\item[]
$\bullet$ the reduction of the even vector field $D_i^2$ gives a
nowhere vanishing (complex) vector field $(D_i^2)_{red}$ on $U_i$.\\
$\bullet$ the complex conjugate of $(D_i^2)_{red}$ is
$-(D_i^2)_{red}$.
\end{itemize}
3) the restrictions of $D_i$ and $D_j$ to $U_i\cap U_j$ are equal up
to a possible sign. \newline Tow such collections define the same
structure if their union is again such a structure.
\end{definition}

Let $S$ be a cs-manifold. Let's define the category RB$^{d|1}_S$,
which is a super analogy to Definition 3.7 of RB$^d_S$.

\begin{definition}{\bf (The category $\mathrm{RB}_S^{d|1}$)}The
objects and morphisms are the followings:

\noindent {\bf objects} of $\mathrm{RB}_S^{d|1}$ are smooth quasi
bundles of cs-manifolds $U\to S$ with fibers of dimension $d|1$
which are equipped with a fiberwise super Riemannian struture.

\noindent {\bf morphisms} from $U$ to $U'$ are embeddings
$U\hookrightarrow U'$ of cs-manifolds that are bundle maps (i.e.
commutes with the projection $S$) and preserve the fiberwise
Riemannian structure.

\end{definition}

Let's present some examples of objects and morphisms.

\noindent {\bf the super point} spt$_S\in$ RB$_S^{1|1}$. The
quadruple
$$ \mathrm{spt}:=(U, Y, U^{+}, U^-)=(S\times \mathbf{R}_{cs}^{1|1}, S\times
\mathbf{R}_{cs}^{0|1}, S\times \mathbf{R}_{cs, -}^{1|1}, S\times
\mathbf{R}_{cs, +}^{1|1})$$ is an object RB$_S^{1|1}$; here
$\mathbf{R}_{cs, \pm}^{1|1}\subset \mathbf{R}_{cs}^{1|1}$ is the
super submanifold whose reduced manifold is
$\mathbf{R}^1_{\pm}\subset \mathbf{R}^1$.

\noindent {\bf the super interval} $I_l^{1|1}\in
\mathrm{RB}^{1|1}_S(\mathrm{spt}, \mathrm{spt})$. For $l\in
\mathbf{R}_{cs, +}^{1|1}(S)$ the pair of bundle maps
$$ \xymatrix@C=0.5cm{
 U=S\times \mathbf{R}_{cs}^{1|1}\ar[rr]^{\ \mathrm{id}} && \Sigma=S\times \mathbf{R}_{cs}^{1|1} &&\ar[ll]_{\ \ l}U=S\times \mathbf{R}_{cs}^{1|1} },$$
is a super Riemannian bordism from spt$_S$ to spt$_S$. We will use
the notation $I_l^{1|1}$ for this morphism. The $l$ in above diagram
is actually $(1\times \mu)\circ (1\times l \times 1)\circ (\Delta
\times 1)$, where $\mu$ is the product for a standard group
structure on $\mathbf{R}^{1|1}_{cs}$:
$$\mu:\mathbf{R}^{1|1}_{cs}\times \mathbf{R}^{1|1}_{cs}\to
\mathbf{R}^{1|1}_{cs}, \ ((t_1, \theta_1), (t_2, \theta_2))\mapsto
(t_1+t_2+\theta_1\theta_2, \theta_1+\theta_2).$$

\begin{definition} {\bf (The super smooth category
$\mathfrak{SRB}^{d|1}$)} The smooth category $\mathfrak{SRB}^{d|1}$
is the functor
$$\mathfrak{SRB}^{d|1}: \mathrm{cs-SMAN}^{op}\rightarrow \mathrm{CAT}$$ which sends a
cs-manifold $S$ to the category $RB_S^{d|1}$ and $f:S\to S'$ a
smooth map between cs-manifolds to the pullback functor $f^*:
RB^{d|1}_S\to RB^{d|1}_{S'}$.
\end{definition}

A {\it vector bundle} over a cs-manifold $S$ is a sheaf of modules
over the structure sheaf $\mathcal{O}_S$ which is locally isomorphic
to the (projective, graded) tensor product $\mathcal{O}_S\otimes V$,
where $V$ is a $\mathbf{Z}/ 2$-graded locally convex vector space.
These modules are equipped with a locally convex topology and the
local isomorphism is bi-continuous. Here $C^\infty(S)$ comes with
its usual Frechet topology. A {\it vector bundle map} is a
continuous map between these sheaves.

One can very similarly define the category TV$_S^{\pm}$ as in
definition 3.9 by using the notion of quasi vector bundles over $S$
and consequently define the super smooth category
$\mathfrak{STV}^{\pm}$.

\begin{definition}{\bf (The super smooth category
$\mathfrak{STV}^{\pm}$)} The super smooth category
$\mathfrak{STV}^{\pm}$ is the functor
$$ \mathfrak{STV}^{\pm}: \mathrm{cs-SMAN}^{op}\to \mathrm{CAT}$$
which sends a cs-manifold $S$ to TV$^{\pm}_S$ and a smooth map
$f':S'\to S$ to the pull-back functor $f^*:TV^{\pm}_S\to
TV^{\pm}_{S'}$.
\end{definition}

\begin{definition} A super symmetric quantum field theory of
dimension $d|1$ is a natural transformation
$$\mathfrak{SF}^{d|1}: \mathfrak{SRB}^{d|1}\to
\mathfrak{STV}^{\pm}$$  such that for any $S\in
\mathrm{cs}-\mathrm{SMAN}$, $\mathfrak{SF}^{d|1}_S$ is a symmetric
monoidal functor, which is compatible with the involution ${}^-$ and
the anti-involution ${}^{\vee}$.
\end{definition}

Let $M$ be a smooth manifold. We can also similarly as what we did
after definition 3.11 define relevant categories relative to $M$ and
define a {\bf super symmetric quantum field theory of dimension
$d|1$ over $M$} to be a natural transformation
$$ \mathfrak{SF}^{d|1}(M): \mathfrak{SRB}^{d|1}(M)\to \mathfrak{STV}^{\pm},$$
such that for any $S\in \mathrm{cs}-\mathrm{SMAN}$,
$\mathfrak{SF}^{d|1}_S(M)$ is a symmetric monoidal functor, which is
compatible with the involution ${}^-$ and the anti-involution
${}^{\vee}$.

In the following two subsections, we talk about some examples of
SUSY QFTs arising from classical geometric objects.
\subsubsection{Super Parallel Transport} Let $E$ be a complex vector
bundle over $M$ and $\nabla^E$ be a connection over $E$. Dumitrescu
has introduced the super parallel transport as follows. His
construction with some modifications can be viewed as an example of
the $1|1$D QFT over $M$ defined above.

Let $c:S\times \mathbf{R}^{1|1}\rightarrow M$ be a family of {\it
supercurves} parametrized by $S$ in $M$. Let $c^*E$ and
$c^*\nabla^E$ be the pull back of the vector bundle and connection
to $S\times \mathbf{R}^{1|1}$ respectively. The vector field
$D=\frac{\partial}{\partial \theta}+\theta \frac{\partial}{\partial
t}$ extends trivially to $S\times \mathbf{R}^{1|1}$. Consider the
derivation $(c^*\nabla^E)_D: \Gamma(c^*E)\rightarrow\Gamma(c^*E).$
An element $\psi$ of $\Gamma(c^*E)$ is called a {\it section of E
along c} and called super parallel if moreover it satisfies
\be(c^*\nabla^E)_D\psi=0.\ee

In local coordinates, one can think of the above equation as a so
called (by Dumitrescu) {\it half-order differential equation}. He
named this because of two reasons: first
$D^2=\frac{\partial}{\partial t}$; secondly, for $2n$ unknowns
functions we only need $n$ values as initial data.

\begin{theorem} (Dumitrescu, [7, Prop. 4.1]) Let $S$ be a supermanifold and $c:S\times
\mathbf{R}^{1|1}\rightarrow M$ be a family of {\it supercurves}
parametrized by $S$ in $M$. Let $\psi_0 \in \Gamma(c^*_{0,0}E)$ be
a section of $E$ along $c_{0,0}:S\rightarrow S\times
\mathbf{R}^{1|1}\rightarrow M$, with the first map the standard
inclusion $i_{0,0}: S\rightarrow S\times \mathbf{R}^{1|1}$. Then
there exists a unique super parallel section $\psi$ of $E$ along
$c$ such that $\psi(0,0)=\psi_0.$
\end{theorem}
This theorem is proved in \cite{D06} by writing the equation (3.3)
in local coordinates and reducing it to a system of first ordinary
differential equations.

Let $S$ be a supermanifold and $(t, \theta)\in
\mathbf{R}_+^{1|1}(S)$ be an $S$-point of $\mathbf{R}_+^{1|1}$.
Consider the super triplet
$$\xymatrix@C=0.5cm{
  S\ar[rr]_{i_{(0,0)}\ \ \ } && S\times \mathbf{R}_+^{1|1} && \ar[ll]^{\ \ \ \ \ \ \ \ \ \ i_{(t,\theta)}\ \ \ } S }$$
with $i_{(0,0)}(s)=(s,0,0)$ and $i_{(t, \theta)}(s)=(s, t(s),
\theta(s))$. Denote this (family of) super intervals by $I_{(t,
\theta)}$. Let $x$ and $y$ be $S$-points of $M$. A super path in $M$
parametrized by $I_{(t, \theta)}$ and with endpoints $x$ and $y$ is
a super curve $c:S\times \mathbf{R}^{1|1}\to M$ with $c\circ
i_{(0,0)}=c(0,0)=x$ and $c\circ i_{(t,\theta)}=c(t, \theta)=y$.

Theorem 3.2 tells us that to any superpath $c:S\times I_{t,
\theta}\rightarrow M$ in $M$, one can associate a bundle map $$
\xymatrix{
 x^*E \ar[rr]^{SP(c)} \ar[dr]
                &  &    y*E \ar[dl]    \\
                & S                 }$$
which is called {\it super parallel tranport}. It can be showed that
the super parallel transport satisfies the following properties (see
details in \cite{D06}):
\newline \noindent {\bf (1)} The correspondence $c\rightarrow SP(c)$ is smooth
and natural in $S$. Smoothness means: if $c$ is a family of smooth
superpaths parametrized by a supermanifold $S$, then the map
$SP(c): c^*_{0,0}E\rightarrow c^*_{t, \theta}E$ is a smooth bundle
map over $S$.
\newline \noindent {\bf (2)} Compatibility under glueing: If $c:I_{t, \theta}\rightarrow
M$ and $c':I_{t', \theta'}\rightarrow M$ are two superpaths in $M$
such that $c'=c\circ R_{t,\theta}$ on some neighborhood $S\times
U$ of $S\times(0,0)\hookrightarrow S\times \mathbf{R}^{1|1}$, with
$U$ an open subsupermanifold in $\mathbf{R}^{1|1}$ containing
$(0,0)$, we have
$$SP(c'\cdot c)=SP(c')\circ SP(c),$$
where $c'\cdot c: I_{t'+t+\theta'\theta,
\theta'+\theta}\rightarrow M$ is obtained from $c$ and $c'$ by
gluing them along their common endpoint. Here $R_{t, \theta}:
S\times \mathbf{R}^{1|1}\rightarrow S\times \mathbf{R}^{1|1}$ is
the right translation by $(t, \theta)$ in the $\mathbf{R}^{1|1}$
direction.
\newline \noindent {\bf (3)} For any superpath
$c:I_{t,\theta}\rightarrow M$, the bundle map
$SP(c):c^*_{0,0}E\rightarrow c^*_{t,\theta}E$ is an isomorphism.
\newline \noindent {\bf (4)} Invariance under geometric reparametrization:
Given $c:I_{t,\theta}\rightarrow M$ a superpath in $M$ and
$\varphi: I_{s,\eta} \rightarrow I_{t,\theta}$ a family of
diffeomorphisms of superintervals that preserve the vertical
distribution, we have
 $ SP(c\circ \varphi)=SP(c).$

We want to point out that some modified version of Florin's super
parallel transport canonically associates a
 $1|1$D QFT to $(E, \nabla^E)$. Let's explain the modifications. Since in SUSY QFT, we
 work with cs-manifolds, the $S$ in Florin's setting should be
 changed to cs-manifold. Also we will have to use $\mathbf{R}_{cs}^{1|1}$
 instead of $\mathbf{R}^{1|1}$. Note that $\mathbf{R}_{cs}^{1|1}$
 has a standard super Riemannian structure, i.e. the odd vector
 field $D_{cs}=\frac{1}{2\pi}\frac{\partial}{\partial
 \theta}-i\theta\frac{\partial}{\partial
 t}$ ($\frac{1}{2\pi}$ is a normalization constant we need in this paper; $\frac{\partial}{\partial
 \theta}-i\theta\frac{\partial}{\partial
 t}$ is the one used in \cite{ST07} and the physics motivation of applying it can be found
 there). Therefore we have to use $D_{cs}$ instead of $D$ is
 Florin's  setting. It's not hard to see that with these modifications,
 Florin's construction goes through and the resulted new super
 parallel transport still satisfies those properties. We will use
 the same notation $SP$ in the following. Another thing we want to
 point out is that to only consider super paths in Florin's setting is
 enough for all super bordisms because super parallel is only a
 local construction and the isometric group preserving the super
 Riemannian structure $\mathrm{Isom}(\mathbf{R}_{cs}^{1|1},
 D_{cs})\cong\mathbf{R}_{cs}^{1|1}$.

 Therefore $(E, \nabla^E)$ canonically gives us a $1|1$D QFT
 $\mathfrak{SF}^{1|1}(M)^{(E, \nabla^E)}$. Let's present some examples to show $\mathfrak{SF}^{1|1}(M)^{(E, \nabla^E)}$
 explicitly. Let $i:S\to U=S\times
 \mathbf{R}_{cs}^{1|1}$ be the standard inclusion $i(s)=(s, 0, 0)$.
 Then $$\mathfrak{SF}^{1|1}_S(M)^{(E, \nabla)}(\alpha: U\to M)=i^*\alpha^* E,$$
 where $\alpha: U\to M$ is a super point over $M$ and $i^*\alpha^*
 E$ is the pull back bundle over $S$. For a super interval
 $(I_{cs, l}^{1|1}, \beta)$ over $M$,
 $$  \xymatrix{
 U=S\times\mathbf{R}_{cs}^{1|1} \ar[r]^{\mathrm{id}}
                &\Sigma= S\times\mathbf{R}_{cs}^{1|1}\ar[d]^{\beta}&  \ar[l]_{l} U=S\times\mathbf{R}_{cs}^{1|1}    \\
                & M                 },$$
$\mathfrak{SF}^{1|1}(M)^{(E, \nabla)}$ sends $(I_{cs, l}^{1|1},
\beta)$ to $$\xymatrix{
  i^*\beta^* E \ar[rr]^{SP((I_{cs, l}^{1|1}, \beta))} \ar[dr]
                &  &    i^*l^* \beta^*E  \ar[dl]  \\
                & S                 }.$$

\subsubsection{$0|1$D Theories} Let $M$ be a smooth manifold.
Hohnhold, Kreck, Stolz and Teichner have also studied {\bf $0|1$D
QFT's over $M$} which in spirit are very similar to the above SUSY
QFT's of dimension $d|1$, $d=1,2$. Let $S$ be a cs-manifold. Let's
first define super categories RB$_S^{0|1}(M)$ and TV$_S^{0}$.

\noindent $\bullet$ RB$_S^{0|1}(M)$: For an cs-manifold $S$, the
objects of RB$_S^{0|1}(M)$ are $S$-families of ``super points in
$M$, i.e. an object of RB$_S^{0|1}(M)$ is a pair $(S\times
\mathbf{R}_{cs}^{0|1}, f)$, where $f:S\times
\mathbf{R}_{cs}^{0|1}\rightarrow M$ is a morphism between
supermanifolds. A morphism from $(S\times \mathbf{R}_{cs}^{0|1}, f)$
to $(S\times \mathbf{R}_{cs}^{0|1}, f')$ is a diffeomorphism $G$
making the following diagram commutative:
\begin{displaymath}
\xymatrix{ & S\times \mathbf{R}_{cs}^{0|1}\ar[dl]_p \ar[dd]^G \ar[dr]^f&\\
          S&                                                      &M\\
           & S\times \mathbf{R}_{cs}^{0|1}\ar[ul]_{p}\ar[ur]_{f'}    &
           }
\end{displaymath}
$$ $$

\noindent $\bullet$ TV$_S^{0}$: For any cs-manifold $S$, the objects
of the category TV$_S^{0}$ are $C^\infty(S)$; the morphisms are
identity morphisms.

Similarly as in the positive dimension cases, we can define a super
smooth category $ \mathfrak{SRB}^{0|1}(M): \mathrm{cs-SMAN} \to
\mathrm{CAT}$ such that $\mathfrak{SRB}^{0|1}(M)(S)=RB_S^{0|1}(M)$
and a super smooth category $\mathfrak{STV}^{0}:\mathrm{cs-SMAN} \to
\mathrm{CAT}$ such that $\mathfrak{STV}^{0}(S)=\mathrm{TV}_S^{0}$.

A {\bf supersymmetric quantum field theory of dimension $0|1$ over
$M$} is a symmetric monoidal functor:
$$\mathfrak{SF}^{0|1}(M): \mathfrak{SRB}^{0|1}(M) \to
\mathfrak{STV}^{0}.$$

The following theorem relates $0|1$D QFT's over $M$ to differential
forms on $M$.
\begin{theorem} (Hohnhold, Kreck, Stolz and Teichner, \cite{HKST06})
There is a bijection
\begin{displaymath} \xymatrix {\{\mathfrak{SF}^{0|1}(M): \mathfrak{SRB}^{0|1}(M)\rightarrow \mathfrak{STV}^{0}\}
\ar[r]^{\ \ \ \ \ \ \ \ e} &\{\omega\in \Omega^*(M)|\,d\omega=0\}},
\end{displaymath} where $e$ is given by
\begin{equation*}
\begin{split}
 &e(\mathfrak{SF}^{0|1}(M))\\
=&\mathfrak{SF}^{0|1}(M)(\underline{\mathbf{SM}}(\mathbf{R}_{cs}^{0|1},
M))(\underline{\mathbf{SM}}(\mathbf{R}_{cs}^{0|1}, M)\times
\mathbf{R}_{cs}^{0|1}, ev)\in
C^\infty(\underline{\mathbf{SM}}(\mathbf{R}_{cs}^{0|1}, M))\cong
\Omega^*(M).\end{split}\end{equation*} Here
$ev:\underline{\mathbf{SM}}(\mathbf{R}_{cs}^{0|1}, M)\times
\mathbf{R}_{cs}^{0|1}\rightarrow M$ is the evaluation map.
\end{theorem}

\begin{remark} There is also an equivariant version of the above
$0|1D$ theory.

\end{remark}

\section{SUSY QFTs and Chern Character} In this section we
construct the Bismut-Chern character form via our modified super
parallel transport and super loop spaces. Our construction gives the
Bismut-Chern character a quantum interpretation. It shows that the
Bismut-Chern character can be viewed as a map from the 1D SUSY QFT
induced by a vector bundle with connection over $M$ to a 0D SUSY QFT
over the loop space of $M$. With the quantum interpretation of the
Bismut-Chern character, we are able to construct Chern character
type maps in the world of SUSY QFTs.

\subsection{The Chern Character and Bismut-Chern Character}
Let $E$ be a vector bundle over a smooth manifold $M$. Let
$\nabla^E$ be a connection on $E$ and $R^E=(\nabla^E)^2$ be the
curvature. The Chern character form associated to $(E, \nabla^E)$ is
defined as
$$ \mathrm{Ch}(E,
\nabla^E):=\mathrm{Tr}\left(e^{\frac{\sqrt{-1}}{2\pi}R^E}\right),$$
which is an even closed differential form on $M$ (by the Chern-Weil
theory). The Chern character form induces the Chern character
homomorphism:
$$\mathrm{Ch}:K(M)\rightarrow H_{dR}^{even}(M, \mathbf{C}).$$ The
importance of this homomorphism lie in the following result due to
Atiyah and Hirzebruch, which says that if one ignores the torsion
elements in $K(M)$, the induced homomorphism:
$$\mathrm{Ch}: K(M)\otimes \mathbf{C}\rightarrow H_{dR}^{even}(M,
\mathbf{C})$$ is actually an isomorphism. The Chern character plays
an important role in geometry and topology.

Let $LM$ be the free loop space of $M$. It is the set of $C^\infty$
mappings $t\in S^1\rightarrow x_t\in M.$ If $x\in M$, the tangent
space $T_xLM$ is identified with the space of smooth periodic vector
fields $X$ over $x$ so that $X_t\in T_{x_t}M$. $LX$ is modeled as a
Frech$\mathrm{\acute{e}}$t manifold.

The Chern character form has a loop space lifting, the Bismut-Chern
character form \cite{B85}. Let's roughly explain the motivation of
this lifting (cf. \cite{GJP91}). Let $S$ be a Clifford module on $M$
with Dirac operator $D$. Witten observed that it should be possible
to associate an equivariantly closed current $\mu_D$ on the free
loop space of $M$ such that $\mathrm{Ind}D=\langle \mu_D, 1\rangle$,
where $1\in \Omega^0(LM)$. The source of this current is the
formalism of path-integrals in supersymmetric quantum mechanics.
Bismut \cite{B85} showed how to generalize this formula by
associating to a vector bundle $E$ over $M$, equipped with a
connection $\nabla^E$, an equivariantly closed differential form
$\mathrm{Bch}(E, \nabla^E)$ on $LM$ such that $\mathrm{Ind}(D\otimes
E)=\langle \mu_D, \mathrm{Bch}(E, \nabla^E)\rangle$, and moreover
$i^*\mathrm{Bch}(E, \nabla^E)=\mathrm{Ch}(E, \nabla^E)$, where
$i:M\to LM$ is the inclusion of the point loops. Let's call
$\mathrm{Bch}(E, \nabla^E)$ the {\it Bismut-Chern character form}
associated to $(E, \nabla^E)$. Getzler, Jones and Petrack
\cite{GJP91} give a formula for the Bismut-Chern character form from
the point of view of their model of equivariant differential forms
on loop space. In their model, they reformulate equivariant
differential forms and equivariant currents on $LM$ as cyclic chains
and cochains over differential graded algebra $\Omega(M)$ of
differential forms on $M$.

The Bismut-Chern character form $\mathrm{BCh}(E, \nabla^E)$ is
defined as follows (\cite{B85}, cf.\cite{GJP91}):
\begin{definition} Let $\varphi(t)\in \Omega(LM)\otimes
\mathrm{Hom}(E_{\gamma(0)}, E_{\gamma(t)})$ be the solution of the
ordinary differential equation \be \nabla_{\frac{\partial}{\partial
t}}\varphi(t)=\varphi(t)\frac{iF(t)}{2\pi},\ee with initial
condition $\varphi(0)=1$, where $F(t)$ is the pull back of the
curvature of $\nabla^E$ at $\gamma(t)$ by the parallel transport
along the loop $\gamma$. Then \be \mathrm{BCh}(E,
\nabla^E):=\mathrm{Tr}(\varphi(1)) \in \Omega(LM).\ee
\end{definition}

Bismut showed that $\mathrm{BCh}(E, \nabla^E)$ defined in this way
is an equivariant closed form over $LM$ and the restriction to $M$,
the fixed point set of $LM$, is the ordinary Chern character form
$\mathrm{Ch}(E, \nabla^E)$ over $M$. In the following we will
provide a new understanding of $\mathrm{BCh}(E, \nabla^E)$ via
modified super parallel transport.

\subsection{Super Loop Spaces}

Let $S_{cs}^{1|1}=S^1\times \mathbf{R}_{cs}^{0|1}$ be the standard
super circle. Let $\underline{\mathbf{SM}}(S_{cs}^{1|1}, M)$ be the
{\it super loop space} of $M$ (which is a generalized super manifold
currently). Let $\Pi TLM_{cs}$ be the supermanifold defined as in
Section 2.1 such that functions on it are sections
$\Lambda^*(T^*LM)\otimes\mathbf{C}$. Since we are now working with
infinite dimensional case, let's explain it a little bit. We equip
the fibre $T_\gamma LM$ with the Frech\'{e}t topology and the fibre
$T_\gamma^*LM$ with the induced weak topology. Note that with this
weak topology, $(T_\gamma^*LM)^*=T_\gamma LM$ (cf. \cite{Con}). Let
$\Lambda^*(T^*LM):=\overline{\oplus_{k=0}^\infty \Lambda^k(T^*LM)}$,
the fibrewise completion of $\oplus_{k=0}^\infty \Lambda^k(T^*LM)$.
We call $\Pi TLM_{cs}$ a complex super Frech\'{e}t manifold.

For any super manifold $S$, by definition, one has \be \begin{split}
&\mathbf{SM}(S,\underline{\mathbf{SM}}(S_{cs}^{1|1}, M))\\
=&\mathbf{SM}(S\times S_{cs}^{1|1}, M)\\
=&\mathbf{SM}(S\times S^1\times \mathbf{R}_{cs}^{0|1}, M)\\
=&\mathbf{SM}(S\times \mathbf{R}_{cs}^{0|1}, LM)\\
=&\{\mathrm{pairs} \ (a, X_a)|a\in LM(S), X_a\in TLM_a\otimes
\mathbf{C}, X_a \
\mathrm{odd}\}\\
=&\mathbf{SM}(S, \Pi TLM_{cs}),
\end{split}\ee
where we can see the last equality holds from the proof of
Proposition 2.1.1 by requiring all the maps between algebras to be
maps between Frech\'{e}t algebras. Therefore we have

\begin{proposition} \be \underline{\mathbf{SM}}(S_{cs}^{1|1}, M)\cong \Pi
TLM_{cs}.\ee
\end{proposition}
This essentially means that $\Pi TLM_{cs}$ is a model for the
generalized super manifold $\mathbf{SM}(- \times S_{cs}^{1|1}, M)$.

\begin{definition} Let $\omega$ be a degree $p$ differential form over $M$.
Given any $t\in S^1$, $\omega$ defines a differential form
$\omega_t$ over $LM$: \be \omega_t|_x(X^1, X^2, \cdots,
X^p):=\omega(X^1_t, \cdots, X^p_t),\ee where $X^1, \cdots, X^p$
are tangent vectors of $LM$ at $x$. If $\omega$ is degree $1$,
given $t\in S^1$, we can also canonically define a smooth function
$\omega(t)$ on $LM$ by \be \omega(t)(x):=\omega(\dot{x}(t)).\ee As
$t$ runs over all of $S^1$, $\omega(t)$ then defines a smooth
function on $LM\times S^1$, which we denote by
$\widetilde{\omega}$.
\end{definition}

As $t$ runs over all of $S^1$, one can view $\omega_t$  defined
above as a $C^\infty$-function on
$\underline{\mathbf{SM}}(S_{cs}^{1|1}, M)\times S^1\cong \Pi
TLM_{cs}\times S^1$. Denote it by $\widehat{\omega}$. Note that
$\widehat{\omega}$ can be viewed as a function on
$\underline{\mathbf{SM}}(S_{cs}^{1|1}, M)\times S_{cs}^{1|1}$.
Similarly, $\widetilde{\omega}$ can also be viewed as a function on
$\underline{\mathbf{SM}}(S_{cs}^{1|1}, M)\times S_{cs}^{1|1}$.

Let $slev:\underline{\mathbf{SM}}(S_{cs}^{1|1}, M)\times
S_{cs}^{1|1}\rightarrow M$ and $lev:LM\times S^{1}\rightarrow M$ be
the evaluation maps defined after (2.3)-(2.5) above. Here $slev$ and
$lev$ represent super loop evaluation and loop evaluation
respectively.

The following theorem characterizes the map $slev$ on function
spaces.
\begin{theorem}The super loop evaluation map $slev:\underline{\mathbf{SM}}(S_{cs}^{1|1}, M)\times
S_{cs}^{1|1}\rightarrow M$ is characterized on functions by
$$C^\infty(M)\rightarrow C^\infty(\underline{\mathbf{SM}}(S_{cs}^{1|1}, M)\times
S_{cs}^{1|1})=C^\infty(\Pi TLM \times S_{cs}^{1|1}),$$ \be f \mapsto
\widehat{f}+\theta (\widehat{df}).\ee
\end{theorem}
\begin{proof}
Let $S$ be any supermanifold, $\varphi\in \mathbf{SM}(S\times
S_{cs}^{1|1}, M), \phi \in \mathbf{SM}(S, S_{cs}^{1|1})$. Let
$$ r:\mathbf{SM}(S\times
S_{cs}^{1|1}, M)\rightarrow \mathbf{SM}(S,
\underline{\mathbf{SM}}(S_{cs}^{1|1}, M))$$ be the representation
map.

Note that the following diagram is commutative:

\begin{displaymath}
\xymatrix{\mathbf{SM}(S\times S_{cs}^{1|1}, M)\times \mathbf{SM}(S,S_{cs}^{1|1})\ar[d]_{r\times id}\ar[r]^{\ \ \ \ \ \ \ \ \ \ \ \ \ slev_S}& \mathbf{SM}(S,M)\\
         \mathbf{SM}(S, \underline{\mathbf{SM}}(S_{cs}^{1|1}, M))\times \mathbf{SM}(S,S_{cs}^{1|1})\ar[r]_{\ \ \ \ \ \  m} & \mathbf{SM}(S,\underline{\mathbf{SM}}(S_{cs}^{1|1}, M)\times S_{cs}^{1|1})\ar[u]_{slev(S)}\\
          }
\end{displaymath}

From (2.4), (2.5), we see that
\begin{displaymath}\xymatrix{S \ar[r]^\Delta &{S\times S}\ar[r]^{1\times\phi} &{S\times S_{cs}^{1|1}}\ar[r]^\varphi &M\ar[r]^f &\mathbf{C}}
\end{displaymath} characterizes $slev_S$, or in other words
\be slev_S^*(f)= f\circ\varphi\circ(1\times \phi)\circ\Delta. \ee

To prove our theorem, we only have to verify that $$
slev_S^*(f)(s)=(\widehat{f}+\theta\widehat{df})(r(\varphi)(s),
\phi(s)),$$ i.e.
$$ f(\varphi(s, \phi(s)))=(\widehat{f}+\theta\widehat{df})(r(\varphi)(s),
\phi(s)).$$

We can actually prove that \be f(\varphi(s,t,
\theta))=(\widehat{f}+\theta\widehat{df})(r(\varphi)(s), t,
\theta)=\widehat{f}(r(\varphi)(s), t)+\theta
\widehat{df}(r(\varphi)(s), t).\ee  This is essentially only an
explicit explanation of the the canonical isomorphism
$$\underline{\mathbf{SM}}(\mathbf{R}_{cs}^{0|1}, M)\cong \Pi
TM_{cs}$$ in our situation.

Let \be f(\varphi(s,t, \theta))=\varphi_1(f)(s,t)+\theta
\varphi_2(f)(s,t),\ee where $\varphi_1(f)(s,t)$ and
$\varphi_2(f)(s,t)$ are two functions on $S\times S^1$. Let $g\in
C^\infty(M)$ be another function. We should have that
$$\varphi_1(fg)(s,t)+\theta \varphi_2(fg)(s,t)=(\varphi_1(f)(s,t)+\theta
\varphi_2(f)(s,t))(\varphi_1(g)(s,t)+\theta \varphi_2(g)(s,t)).$$
Therefore
$$\varphi_1(fg)=\varphi_1(f)\varphi_1(g),$$
$$\varphi_2(fg)=\varphi_1(f)\varphi_2(g)+\varphi_1(g)\varphi_2(f).$$

We can see from above that $\varphi$ gives two maps $\varphi_1,
\varphi_2:C^\infty(M)\rightarrow C^\infty(S\times S^1)$ such that
$\varphi_1$ is a homomorphism of rings and $\varphi_2$ is
$\varphi_1\circ X_\varphi$, where $X_\varphi$ is a tangent vector
field over $M$.

On the other hand, by the definition of the representation map:
$$ r:\mathbf{SM}(S\times
S_{cs}^{1|1}, M)\rightarrow \mathbf{SM}(S,
\underline{\mathbf{SM}}(S_{cs}^{1|1}, M)),$$ it's not hard to see
that
$$\widehat{f}(r(\varphi)(s),
t)=f(\varphi(s,t,0))=\varphi_1(f)(s,t),$$ and
$$\widehat{df}(r(\varphi)(s), t)=X_\varphi(f)(\varphi(s,t,0))=\varphi_1(X_\varphi(f))(s,t)=\varphi_2(f)(s,t).$$

This finishes the proof.
\end{proof}
We have the following diagram of maps (which does not commute):
\begin{displaymath}
\xymatrix{{\underline{\mathbf{SM}}(S_{cs}^{1|1}, M)\times S_{cs}^{1|1}}\ar[d]_{p}\ar[dr]^{slev}& \\
           LM\times S^{1} \ar[r]_{lev}                                        &M }
\end{displaymath}

Let's adopt the Deligne-Morgan sign convention in the following. The
paring of vectors and $1$-forms is written with the vector on the
left with the rule
$$\langle  uD, v\omega \rangle=(-1)^{p(D)p(v)}uv\langle D, \omega\rangle.
$$ When computing with differential forms on super space, we sue a
bigraded point of view. This means objects have a ``cohomological
degree" and a parity. The permutation of objects of parity $p_1,
p_2$ and cohomological degree $d_1, d_2$ introduce two signs
$(-1)^{d_1d_2}$ and an additional factor $(-1)^{p_1p_2}$. See
\cite{DM99} for details.

We find that $slev$ also has the following property on differential
one forms.
\begin{theorem} Let $\omega\in \Omega^1(M)$, then one has
\be \langle D_{cs},
slev^*(\omega)\rangle=\frac{1}{2\pi}\widehat{\omega}+\frac{1}{2\pi}\theta(\widehat{d\omega})-i\theta\widetilde{\omega}\in
C^\infty(\underline{\mathbf{SM}}(S_{cs}^{1|1}, M)\times
S_{cs}^{1|1}),\ee where
$D_{cs}=\frac{1}{2\pi}\frac{\partial}{\partial \theta}-i\theta
\frac{\partial}{\partial t}.$
\end{theorem}
\begin{proof} By partition of unity, $\omega$ can be written as
$$\omega=\sum_{i}f_idg_i,$$
where $f_i$'s, $g_i$'s are smooth functions over $M$, $f_i$'s are
nonnegative, $\sum_i f_i=1$ and near each point, there exists a
neighborhood such that there are only finite many $f_i$'s nonzero
in it. Therefore we only have to prove the theorem for
differential forms like $fdg$.

Actually we have \be \begin{split}\langle D_{cs},
slev^*(fdg)\rangle=&\left\langle\frac{1}{2\pi}\frac{\partial}{\partial
\theta}-i\theta \frac{\partial}{\partial t}, (\widehat{f}+\theta
\widehat{df})d(\widehat{g}+\theta \widehat{dg})
\right\rangle\\
=&\frac{1}{2\pi}\left\langle \frac{\partial}{\partial \theta},
(\widehat{f}+\theta \widehat{df})d(\widehat{g}+\theta
\widehat{dg})\right\rangle-i\left\langle\theta
\frac{\partial}{\partial t}, (\widehat{f}+\theta
\widehat{df})d(\widehat{g}+\theta \widehat{dg})\right\rangle\\
=&\frac{1}{2\pi}\widehat{fdg}+\frac{1}{2\pi}\theta\widehat{dfdg}-i\theta \left\langle \frac{\partial}{\partial t}, \widehat{f}d\widehat{g}\right\rangle\\
=&\frac{1}{2\pi}\widehat{fdg}+\frac{1}{2\pi}\theta
\widehat{d(fdg)}-i\theta
\widehat{f}\left(\frac{\partial \widehat{g}}{\partial t}\right)\\
=&\frac{1}{2\pi}\widehat{fdg}+\frac{1}{2\pi}\theta
\widehat{d(fdg)}-i\theta \widetilde{fdg}.
\end{split}\ee
This completes the proof.
\end{proof}

\subsection{A Quantum Interpretation of the Bismut-Chern Character and Chern Character} Applying the super parallel transport equation (3.3),
one can identify the ordinary differential equation that the Bosonic
part of the super parallel transport of $slev$ should satisfy in
local coordinates.
\begin{theorem} Locally suppose $\nabla^E=d+A$, where $A\in \Omega^1(M, \mathrm{End}(E))$ and let $R=(d+A)^2$ be the curvature.
Then the Bosonic super parallel transport $SP(slev, t, 0)$ satisfies
the following ordinary differential equation \be \frac{d}{d
t}SP(slev, t, 0)=-\frac{i}{2\pi}\widehat{R}SP(slev, t, 0)-
\widetilde{A}SP(slev, t, 0).\ee
\end{theorem}
\begin{proof} By definition of super parallel transport, $$(slev^*\nabla^E)_{D_{cs}}SP(slev, t,
\theta)=0.$$ We therefore have \be \left\langle
\frac{1}{2\pi}\frac{\partial}{\partial \theta}-i\theta
\frac{\partial}{\partial t},  dSP(slev, t, \theta)+(slev^*A)SP(slev,
t, \theta)\right\rangle=0.\ee

Let $SP(slev, t, \theta)=SP(slev, t, 0)+\theta SP'(slev, t)$. Then
we have \be\left\langle \frac{1}{2\pi}\frac{\partial}{\partial
\theta}-i\theta \frac{\partial}{\partial t}, dSP(slev, t,
\theta)\right\rangle=\frac{1}{2\pi}SP'(slev, t)-i\theta \frac{d}{d
t}SP(slev, t, 0),\ee and by Theorem 4.2 \be
\begin{split}
&\left\langle\frac{1}{2\pi}\frac{\partial}{\partial
\theta}-i\theta \frac{\partial}{\partial t}, (slev^*A)SP(slev, t, \theta)\right\rangle \\
=&\left(\frac{1}{2\pi}\widehat{A}+\frac{1}{2\pi}\theta\widehat{dA}-i\theta
\widetilde{A}\right)(SP(slev, t,
0)+\theta SP'(slev, t))\\
=&\frac{1}{2\pi}\widehat{A}SP(slev, t,
0)+\theta\left[-\frac{1}{2\pi}\widehat{A}SP'(slev,
t)+\left(\frac{1}{2\pi}\widehat{dA}-i\widetilde{A}\right)SP(slev, t,
0)\right].
\end{split} \ee

From (4.14)-(4.16), we see that $$SP'(slev,
t)=-\widehat{A}SP(slev, t, 0),$$
$$\frac{d}{d t}SP(slev, t, 0)=\frac{i}{2\pi}\widehat{A}SP'(slev,
t)-\left(\frac{i}{2\pi}\widehat{dA}+\widetilde{A}\right)SP(slev, t,
0).$$

Hence we obtain that
$$\frac{d}{d t}SP(slev, t, 0)=-\frac{i}{2\pi}(\widehat{A^2}+\widehat{dA})SP(slev, t, 0)-\widetilde{A}SP(slev, t, 0), $$
or \be\frac{d}{d t}SP(slev, t, 0)=-\frac{i}{2\pi}\widehat{R}SP(slev,
t, 0)-\widetilde{A}SP(slev, t, 0).\ee
\end{proof}

From Theorem 4.3, we see that super parallel transport along the
fermionic direction $D_{cs}$ on the super evaluation curve is
different from the ordinary transport along the time direction.
However we will see that the super parallel transport along the
curve $lev \circ p$ degenerates to the ordinary transport on
$LM\times S^1$.
\begin{theorem} Using the same notations as in theorem 4.3, we
have $$ \frac{d}{d t}SP(lev\circ p, t, 0)+\widetilde{A}SP(lev\circ
p, t, 0)=0.$$
\end{theorem}
\begin{proof}
It's not hard to see that the map $lev\circ p:
\underline{\mathbf{SM}}(S_{cs}^{1|1}, M)\times
S_{cs}^{1|1}\rightarrow M$ is characterized on functions by
$$C^\infty(M)\rightarrow C^\infty(\underline{\mathbf{SM}}(S_{cs}^{1|1}, M)\times
S_{cs}^{1|1}),  \ \ \ \ \ \ f \mapsto \widehat{f}.$$ Similarly to
what we did in Theorem 4.2, we have

\begin{equation*} \begin{split}\langle D_{cs},
(lev\circ
p)^*(fdg)\rangle=&\left\langle\frac{1}{2\pi}\frac{\partial}{\partial
\theta}-i\theta \frac{\partial}{\partial t}, \hat{f}d\hat{g}
\right\rangle\\
=&-i\theta \left\langle \frac{\partial}{\partial t}, \widehat{f}d\widehat{g}\right\rangle\\
=&-i\theta
\widehat{f}\left(\frac{\partial \widehat{g}}{\partial t}\right)\\
=&-i\theta \widetilde{fdg}.
\end{split}\end{equation*}
Therefore one has for $\omega\in \Omega^1(M)$, $$ \langle D_{cs},
(lev\circ p)^*(\omega)\rangle=-i\theta\widetilde{\omega}\in
C^\infty(\underline{\mathbf{SM}}(S_{cs}^{1|1}, M)\times
S_{cs}^{1|1}).$$

By definition of super parallel transport along the curve $lev\circ
p$,
$$((lev\circ p)^*\nabla^E)_{D_{cs}}SP(lev\circ p, t, \theta)=0.$$ We therefore
have $$ \left\langle \frac{1}{2\pi}\frac{\partial}{\partial
\theta}-i\theta \frac{\partial}{\partial t},  dSP(lev\circ p, t,
\theta)+((lev\circ p)^*A)SP(lev\circ p, t, \theta)\right\rangle=0.$$
Let $SP(lev\circ p, t, \theta)=SP(lev\circ p, t, 0)+\theta
SP'(lev\circ p, t)$. Then we have $$ \left\langle
\frac{1}{2\pi}\frac{\partial}{\partial \theta}-i\theta
\frac{\partial}{\partial t}, dSP(lev\circ p, t,
\theta)\right\rangle=\frac{1}{2\pi}SP'(lev\circ p, t)-i\theta
\frac{d}{d t}SP(lev\circ p, t, 0),$$ and \begin{equation*}
\begin{split}
&\left\langle\frac{1}{2\pi}\frac{\partial}{\partial
\theta}-i\theta \frac{\partial}{\partial t}, ((lev\circ p)^*A)SP((lev\circ p), t, \theta)\right\rangle \\
=&\left(-i\theta \widetilde{A}\right)(SP(lev\circ p, t,
0)+\theta SP'(lev\circ p, t))\\
=&\left(-i\theta \widetilde{A}\right)SP(lev\circ p, t, 0).
\end{split} \end{equation*}

Hence we see that $$SP'(lev\circ p, t)=0,$$
$$\frac{d}{d t}SP(lev\circ p, t, 0)+\widetilde{A}SP(lev\circ p, t,
0)=0.$$
\end{proof}

Combining the above super parallel transports along the two super
curves, we have

\begin{theorem} Let's use the same notations as in theorem 4.3.
The following ordinary differential equations hold: \be
\begin{split} &\frac{d}{d t}[SP^{-1}(lev\circ p, t, 0)SP(slev, t,
0)]\\
=&-SP^{-1}(lev\circ p, t,
0)\left(\frac{i}{2\pi}\widehat{R}\right)[SP^{-1}(lev\circ p, t,
0)SP(slev, t, 0)],\end{split}\ee \be
\begin{split} &\frac{d}{d t}[SP^{-1}(slev, t,
0)SP(lev\circ p, t, 0)]\\
=&[SP^{-1}(slev, t, 0)SP(lev\circ p, t, 0)]SP^{-1}(lev\circ p, t,
0)\left(\frac{i}{2\pi}\widehat{R}\right).\end{split}\ee

\end{theorem}
\begin{proof} Differentiating the identity $SP^{-1}(lev\circ p, t, 0)SP(lev\circ p, t,
0)=I$, we have that
$$\left( \frac{d}{d t}SP^{-1}(lev\circ p, t, 0)\right)SP(lev\circ p, t,
0)=-SP^{-1}(lev\circ p, t, 0)\left(\frac{d}{d t}SP(lev\circ p, t,
0)\right).$$ However by Theorem 4.4, we see that
$$ \frac{d}{d t}SP(lev\circ p,
t, 0)=-\widetilde{A}SP(lev\circ p, t, 0). $$ Therefore, one has \be
\frac{d}{d t}SP^{-1}(lev\circ p, t, 0)= SP^{-1}(lev\circ p, t,
0)\widetilde{A}.\ee

Hence by (4.13) and (4.20), we obtain that \be \begin{split}
&\frac{d}{d t}[SP^{-1}(lev\circ p, t,
0)SP(slev, t, 0)]\\
=&\left(\frac{d}{d t}SP^{-1}(lev\circ p, t, 0)\right)SP(slev, t,
0)+SP^{-1}(lev\circ p, t, 0)\frac{d}{d t}SP(slev, t, 0)\\
=&SP^{-1}(lev\circ p, t, 0)\widetilde{A}SP(slev, t, 0)
\\&+SP^{-1}(lev\circ p, t, 0)\left[-\frac{i}{2\pi}\widehat{R}SP(slev, t, 0)-
\widetilde{A}SP(slev, t, 0)\right]\\
=&-SP^{-1}(lev\circ p, t,
0)\circ\left(\frac{i}{2\pi}\widehat{R}\right)\circ SP(slev, t,
0)\\
=&-\left[SP^{-1}(lev\circ p, t,
0)\circ\frac{i}{2\pi}\widehat{R}\circ SP(lev\circ p, t,
0)\right][SP^{-1}(lev\circ p, t, 0)SP(slev,
t, 0)]\\
=&-SP^{-1}(lev\circ p, t,
0)\left(\frac{i}{2\pi}\widehat{R}\right)[SP^{-1}(lev\circ p, t,
0)SP(slev, t, 0)],
\end{split}\ee which proves (4.18).

From (4.18), it's not hard to obtain (4.19) as how we obtained
(4.20).
\end{proof}

Comparing Definition 4.1 with Theorem 4.5, we obtain that
\begin{theorem} The following identity holds:
\be \mathrm{BCh}(E, \nabla^E)=\mathrm{Tr}[SP^{-1}(slev, 1, 0)\circ
SP(lev\circ p, 1, 0)]\in C^\infty(\underline{\mathbf{SM}}(S^{1|1},
M))=\Omega^*(LM).\ee Up to some sign, the following identity holds:
$$ \mathrm{BCh}(E, \nabla^E)=\mathrm{Tr}[SP^{-1}(lev\circ p, 1, 0)\circ
SP(slev, 1, 0)]\in C^\infty(\underline{\mathbf{SM}}(S_{cs}^{1|1},
M))=\Omega^*(LM).$$
\end{theorem}

\begin{remark} In certain sense, $SP^{-1}(slev, 1,
0)\circ SP(lev\circ p, 1, 0)$ is a loop-deloop process. Therefore,
from Theorem 4.6, we see that the Bismut-Chern character is a
phenomena related to this loop-deloop process when one moves from
$1|1$D theories down to $0|1$D theories. This theorem also shows us
the supersymmetric aspect of the Bismut-Chern character form.
\end{remark}

Let $i:\underline{\mathbf{SM}}(\mathbf{R}_{cs}^{0|1}, M)\to
\underline{\mathbf{SM}}(S_{cs}^{1|1}, M)$ be the inclusion of super
constant loops. Let $$\xymatrix@C=0.5cm{
 \underline{\mathbf{SM}}(\mathbf{R}_{cs}^{0|1}, M)\times S_{cs}^{1|1} \ar[r]^{i\times 1}
 & \underline{\mathbf{SM}}(S_{cs}^{1|1}, M)\times S_{cs}^{1|1} \ar[r]^{\ \ \ \ \ \ \ \ \ \ \ \ slev} & M }$$
be the restriction of the super evaluation curve on the constant
loops. Then it's not hard to see that
\begin{theorem}The following identity holds:
\be \mathrm{Ch}(E, \nabla^E)=\mathrm{Tr}[SP^{-1}(slev\circ (i\times
1), 1, 0)]\in
C^\infty(\underline{\mathbf{SM}}(\mathbf{R}_{cs}^{0|1},
M))=\Omega^*(M).\ee Up to some sign, the following identity holds:
$$ \mathrm{Ch}(E, \nabla^E)=\mathrm{Tr}[SP(slev\circ (i\times
1), 1, 0)]\in
C^\infty(\underline{\mathbf{SM}}(\mathbf{R}_{cs}^{0|1},
M))=\Omega^*(M).$$
\end{theorem}

\subsection{Chern Character in SUSY QFTs} In view of Theorem 4.6, Remark 4.1 and Theorem 4.7, we are motivated
to define Chern character type maps for SUSY QFTs in the following.

Given any $1|1D$ QFT $\mathfrak{SF}^{1|1}(M)$ over $M$, let's
construct a $0|1$D QFT $\mathfrak{SF}^{0|1}(LM)$ over $LM$ as
follows.

Let $\widetilde{r}: LM\to LM$ be the reverse map, which sends
$\gamma: S^1\to M$ to $$\xymatrix@C=0.5cm{
  S^1 \ar[r]^{r} & S^1 \ar[r]^\gamma & M }, $$ where $r(t)=1-t$.
Then $\widetilde{r}$ induces a map, we still denote it by
$\widetilde{r}$, on the super loop space $\widetilde{r}:
\underline{\mathbf{SM}}(\mathbf{R}_{cs}^{0|1},LM)\to
\underline{\mathbf{SM}}(\mathbf{R}_{cs}^{0|1},LM)$. Let's pick out
two particular bordisms:

$$  \xymatrix{
 \underline{\mathbf{SM}}(\mathbf{R}_{cs}^{0|1},LM)\times\mathbf{R}_{cs}^{1|1} \ar[r]^{\mathrm{id}}
                &\underline{\mathbf{SM}}(\mathbf{R}_{cs}^{0|1},LM)\times\mathbf{R}_{cs}^{1|1}\ar[d]^{\widetilde{r}\times \mathrm{id}}&  \ar[l]_{1} \underline{\mathbf{SM}}(\mathbf{R}_{cs}^{0|1},LM)\times\mathbf{R}_{cs}^{1|1}    \\
                & \underline{\mathbf{SM}}(\mathbf{R}_{cs}^{0|1},LM)\times\mathbf{R}_{cs}^{1|1}\ar[d]^{slev}                 \\
                &M},$$
and $$  \xymatrix{
 \underline{\mathbf{SM}}(\mathbf{R}_{cs}^{0|1},LM)\times\mathbf{R}_{cs}^{1|1} \ar[r]^{\mathrm{id}}
                &\underline{\mathbf{SM}}(\mathbf{R}_{cs}^{0|1},LM)\times\mathbf{R}_{cs}^{1|1}\ar[d]^{p}&  \ar[l]_{1} \underline{\mathbf{SM}}(\mathbf{R}_{cs}^{0|1},LM)\times\mathbf{R}_{cs}^{1|1}    \\
                & LM\times S^1\ar[d]^{lev}                 \\
                &M},$$
where
$1:\underline{\mathbf{SM}}(\mathbf{R}_{cs}^{0|1},LM)\times\mathbf{R}_{cs}^{1|1}
\to
\underline{\mathbf{SM}}(\mathbf{R}_{cs}^{0|1},LM)\times\mathbf{R}_{cs}^{1|1}
$ is the constant map given by $1(s, t, \theta)=(s, t+1, \theta).$
Let's denote these two bordisms by $b(1, \widetilde{r}, slev)$ and
$b(1, p, lev)$ respectively.

Let $S$ be any cs-manifold. An object $ f\in \mathbf{SM}(S\times
\mathbf{R}_{cs}^{0|1}, LM)$ determines a map $\widetilde{f} \in
\mathbf{SM}(S,
\underline{\mathbf{SM}}(\mathbf{R}_{cs}^{0|1},LM))=\mathbf{SM}(S,
\underline{\mathbf{SM}}(S_{cs}^{1|1},M))$. With this
$\widetilde{f}$, one has two new bordisms (of $S$-families):

$$  \xymatrix{S\times\mathbf{R}_{cs}^{1|1}\ar[d]_{\widetilde{f}\times \mathrm{id}}\ar[r]^{\mathrm{id}} & S\times\mathbf{R}_{cs}^{1|1} \ar[d]^{\widetilde{f}\times \mathrm{id}}&S\times\mathbf{R}_{cs}^{1|1} \ar[d]^{\widetilde{f}\times \mathrm{id}}\ar[l]_{\mathrm{1}} \\
 \underline{\mathbf{SM}}(\mathbf{R}_{cs}^{0|1},LM)\times\mathbf{R}_{cs}^{1|1} \ar[r]^{\mathrm{id}}
                &\underline{\mathbf{SM}}(\mathbf{R}_{cs}^{0|1},LM)\times\mathbf{R}_{cs}^{1|1}\ar[d]^{\widetilde{r}\times \mathrm{id}}&  \ar[l]_{1} \underline{\mathbf{SM}}(\mathbf{R}_{cs}^{0|1},LM)\times\mathbf{R}_{cs}^{1|1}    \\
                & \underline{\mathbf{SM}}(\mathbf{R}_{cs}^{0|1},LM)\times\mathbf{R}_{cs}^{1|1}\ar[d]^{slev}                 \\
                &M}$$ and
$$  \xymatrix{S\times\mathbf{R}_{cs}^{1|1}\ar[d]_{\widetilde{f}\times \mathrm{id}}\ar[r]^{\mathrm{id}} &S\times\mathbf{R}_{cs}^{1|1}\ar[d]^{\widetilde{f}\times \mathrm{id}} &S\times\mathbf{R}_{cs}^{1|1} \ar[d]^{\widetilde{f}\times \mathrm{id}} \ar[l]_{\mathrm{1}} \\
 \underline{\mathbf{SM}}(\mathbf{R}_{cs}^{0|1},LM)\times\mathbf{R}_{cs}^{1|1} \ar[r]^{\mathrm{id}}
                &\underline{\mathbf{SM}}(\mathbf{R}_{cs}^{0|1},LM)\times\mathbf{R}_{cs}^{1|1}\ar[d]^{p}&  \ar[l]_{1} \underline{\mathbf{SM}}(\mathbf{R}_{cs}^{0|1},LM)\times\mathbf{R}_{cs}^{1|1}    \\
                & LM\times S^1\ar[d]^{lev}                 \\
                &M}$$  Denote them by $b(\widetilde{f}, 1, \widetilde{r}, slev)$ and
$b(\widetilde{f}, 1, p, lev)$ respectively.

Define \be \begin{split} &\mathfrak{SF}^{0|1}_{S}(LM)((S\times
\mathbf{R}_{cs}^{0|1},
f))\\
:=&\mathrm{Str}\left[\mathfrak{SF}^{1|1}_S(M)(b(\widetilde{f}, 1,
\widetilde{r}, slev))\circ \mathfrak{SF}^{1|1}_S(M)(b(\widetilde{f},
1, p, lev))\right]\in {C^\infty(S)}.\end{split}\ee It's not hard to
check that this indeed gives us a $0|1$D QFT over $LM$. In other
words, we have canonically constructed a loop-deloop map:

$$
LD: \xymatrix{\{\mathrm {1|1D\  QFTs \ over} \ M\}\ar[r]&
\{\mathrm{\ 0|1D \ QFTs \ over} \ LM \}},
$$
which makes the following diagram commutative:

$$
\xymatrix{\{\mathrm {1|1D\ QFTs \ over} \ M\}\ar[r]^{LD}& \{\mathrm {0|1D\  QFTs \ over} \ LM \}\ar[d]^{HKST}\\
         \{\mathrm{vector\  bundles \ with \ connections \ over}\ M \}\ar[u]^{SP}\ar[dr]_{\mathrm{Ch}}\ar[r]_{\ \ \ \ \ \ \ \ \ \mathrm{BCh}} & \{S^1-\mathrm{closed\ forms\  on}\  LM\}\ar[d]^{res}\\
                                                                                                     & \{\mathrm{closed\ forms\  on}\  M\}
         }
$$

Let $i:M\to LM$ be the inclusion of constant loops and $i^*$ be the
pull back of field theories on $LM$ to field theories on $M$. Then
the map $i^* \circ LD$ from $1|1$D SUSY QFTs to $0|1$D SUSY QFTs
plays the role of the Chern character in the framework of SUSY QFTs.

\section{Acknowledgement} The paper was finished when the author was visiting
the Max-Planck Instit$\mathrm{\ddot{u}}$t f$\mathrm{\ddot{u}}$r
Mathematik at Bonn. The author is very grateful to Professor Peter
Teichner for his guidance and helps. He also thanks Professor
Stephan Stolz for many helpful discussions.

\end{document}